%% file: main.tex
\pdfoutput=1
\documentclass{article}

\usepackage[nonatbib, preprint]{neurips_2021}

\usepackage[utf8]{inputenc} 
\usepackage[T1]{fontenc}    
\usepackage{hyperref}       
\usepackage{url}            
\usepackage{booktabs}       
\usepackage{amsfonts}       
\usepackage{nicefrac}       
\usepackage{microtype}      
\usepackage{xcolor}         
\usepackage{amsmath}
\usepackage{cleveref}
\usepackage{graphicx}
\usepackage{float}
\usepackage{subfigure}
\usepackage{multirow}
\usepackage{algorithm}
\usepackage{algorithmic}
\usepackage{caption}
\usepackage{subfigure}

\newcommand{\Yordle}{\texttt{YORDLE}}
\title{\Yordle: An Efficient Imitation Learning for \\ Branch and Bound}

\author{%
  Qingyu Qu \\
  Equal contributions \\
  Beihang University \\
  \texttt{quqingyu@buaa.edu.cn} \\
  \And
   Xijun Li \\
   Equal contributions \\
   MIRA Lab, USTC \\
   \texttt{xijunli@miralab.ai} \\
   \And
   Yunfan Zhou \\
   Chinese University of Hong Kong (Shenzhen) \\
   \texttt{yunfanzhou@link.cuhk.edu.cn} \\
}
\begin{document}

\maketitle
\begin{abstract}
Combinatorial optimization problems have aroused extensive research interests due to its huge application potential. In practice, there are highly redundant patterns and characteristics during solving the combinatorial optimization problem, which can be captured by machine learning models. Thus, the 2021 NeurIPS Machine Learning for Combinatorial Optimization (ML4CO) competition is proposed with the goal of improving state-of-the-art combinatorial optimization solvers by replacing key heuristic components with machine learning
techniques.
This work presents our solution and insights gained by team qqy in the dual task of the competition. Our solution is a highly efficient imitation learning framework for performance improvement of Branch and Bound (B\&B), named \Yordle. It employs a hybrid sampling method and an efficient data selection method, which not only accelerates the model training but also improves the decision quality during branching variable selection. 
In our experiments, \Yordle~ greatly outperforms the baseline algorithm adopted by the competition while requiring significantly less time and amounts of data to train the decision model. Specifically, we use only $1/4$ of the amount of data compared to that required for the baseline algorithm, to achieve around $50\%$ higher score than baseline algorithm. The proposed framework \Yordle~ won the championship of the student leaderboard.
\end{abstract}

\input{content/introduction}
\input{content/background}
\input{content/solution}
\input{content/evaluation}
\input{content/conclusion}

\bibliographystyle{unsrt}
\bibliography{reference}

\end{document}

%% file: content/introduction.tex
\section{Introduction}

In recent years, the combinatorial optimization problem (COP) has aroused extensive interests from both industry and academia since many real-life critical optimization problems such as scheduling, planning, bin packing, etc. can be formulated as COP. People can benefit a lot from solving these problems. 
In general, the COP aims to find an optimal configure in discrete spaces, and is one of the most common mathematical topics in industrial applications. 
Many traditional algorithms are proposed to solve kinds of COPs~\cite{achterberg2005branching}. Through many years of practice, kinds of COPs and corresponding algorithms have been verified their effectiveness in optimizing real-life problems. In practical scenarios, there are many highly similar COPs which are only slightly different in coefficients. For example, managing a large-scale production planning requires solving very similar COPs on a daily basis, with a fixed logistic network and bill of material (BOM) while only the demand changes over time. This change of demand is hard to capture by hand-engineered expert rules, and ML-enhanced approaches offer a possible solution to detect typical patterns in the demand history, to further enhance traditional combinatorial optimization algorithms. Therefore, there is a trend of using machine learning (ML) techniques to boost above algorithms~\cite{khalil2016learning} recently.


Among many kinds of COPs, we mainly consider the Mixed Integer Linear Programming (MILP) since its wide adoption in practical modeling. The branch-and-bound (B\&B) algorithm is placed at the core of solving MILP.
B\&B algorithm enumerates the candidate solutions systematically by means of state space search, where the set of candidate solutions is considered to form a search tree with the full set at the root.
The efficiency of B\&B algorithm mainly depends on branching variable selection and node selection. 
Usually, choosing good variables to branch on can lead to a dramatic reduction in terms of the number of nodes needed to solve an instance.
At present, there is still no universally accepted method for the strategy of branching variable selection.
Traditional methods are mostly simple heuristic rules, such as the Most infeasible branching, Pseudocost branching (PC), Strong Branching (SB), Hybrid Strong/Pseudocost branching, Pseudocost branching with strong branching initialization, Reliability branching, etc.~\cite{achterberg2005branching}.
Alvarez et al. adopted machine learning algorithms early to learn the strategies of branching variable selection in the B\&B algorithm~\cite{Alvarez14asupervised}.
Such kind of learning-based strategies are also known as learning to branch.
Learning to branch is different from the traditional optimization methods.
It introduces the concept of learning in the optimization process to help search the optimal solution more effectively.
Balcan has shown empirically and theoretically that it is possible to learn high-performing branching strategies for a given application domain~\cite{Balcan:5}.
Learning branching policies for MILP has become an active research area.
Most relevant researches use supervised or imitation learning to imitate SB method and specialize it to distinct classes of problems.

Supervised machine learning and imitation learning are currently the mainstream approaches for learning to branch.
Alvarez et al. proposed a new approach that uses supervised learning to improve the performances of optimization algorithms in the context of MILP~\cite{Alvarez14asupervised}.
Khalil et al. proposed a machine learning framework for variable branching in MILP. Based on the data collected by SB method, they learned an easy-to-evaluate surrogate function that mimics the SB method, by means of solving a learning-to-rank problem~\cite{khalil2016learning}. And it is competitive with a state-of-the-art commercial solver.
Gasse et al. proposed a new graph convolutional neural network (GCNN) model for learning to branch, which leverages the natural variable-constraint bipartite graph representation of MILP. They trained the GCNN model via imitation learning from the SB method, and demonstrated that this model produced policies that improved upon state-of-the-art machine learning methods for branching~\cite{gasse2019exact}.
Gupta et al. proposed a new hybrid architecture for efficient branching on CPU machines, which combined the expressive power of GCNNs with computationally inexpensive multi-layer perceptrons (MLPs) for branching~\cite{gupta2020hybrid}.
More related researches can refer to the survey provided by Huang et al~\cite{huang2021branch}.
However, \textit{these algorithms often require a large amount of expert data. Considering that there is no commonly accepted expert rules for B\&B, improper expert data may be misleading for training.}
To address this issue, we propose a learning-based framework to generate high-performance expert data for a given metric.

In this work, we develop a highly efficient imitation learning framework for B\&B, named \Yordle~\footnote{\Yordle~ are a race of spirits in a game named `League of Legends' developed by Riot Games. Usually, \Yordle~ are much smaller than humans. Furthermore, they are commonly characterized by their creation and utilization of complex tools. In a word, \Yordle~ are usually small but powerful, which is consistent with the characteristics of our framework.}, which not only greatly accelerates the learning convergence but also improve the quality of branching variable selection. Specifically, the framework consists of four parts, namely data collection, data selection, model learning, and ML-based branching.
In the data collection phase, we adopt a hybrid strategy to collect expert data. And then in the data selection phase, we select the state-action pairs that possess higher cumulative reward. In the model learning phase, a graph convolutional neural network is adopted to train a branching strategy. Finally in the ML-based branching phase, the trained branching strategy is evaluated in parallel.
The technical contributions are summarized as follows:
\begin{itemize}
    \item A hybrid sampling method based on pseudo cost and active constraint method is designed to collect data for demonstration. 
    The demonstration data will be taken as input into a graph neural network (GNN) to train a branching policy.
    \item Before training branching policy, the demonstration data is filtered via a Best-Action Imitation Learning algorithm (BAIL) which selects state-action pairs that possess higher cumulative reward among all state-action pairs within the demonstration data, which results in faster learning convergence.
    \item \Yordle~ is tested over practical MILP dataset obtained from Combinatorial Optimization (ML4CO) NeurIPS 2021 competition~\cite{ml4co}, which greatly outperforms the state-of-the-art imitation learning-based B\&B algorithm with respect to dual integral. 
    \item In ML4CO competition, we only submit two branching model (for `Item Placement' and `Anonymous' dataset respectively) trained using \Yordle~ and the branching model for `Load Balancing' dataset still adopts the baseline algorithms, which won the championship of the student leaderboard in the dual task.
    \item After the competition, we continue to test our framework over `Load Balancing' dataset, which suggests highly competitive performance compared to winner of global leaderboard.
\end{itemize}

Beside, two key insights are gained during the competition, which are as following:

\noindent\textbf{Remark 1.} It is commonly believed that the strong branching strategy can lead to the smallest B\&B tree. However, when the dual integral is adopted as the metric, the data collected by our proposed framework leads to better result than that collected by strong branching strategy. Therefore, we doubt whether the strong branching strategy is still the `golden standard'. 

\noindent\textbf{Remark 2.} During the training process, it is found that a smaller cross-entropy loss does not necessarily lead to a better score, which is worth studying in the future.

%% file: content/background.tex
\section{Background}

\subsection{Mixed integer linear programs}

A mixed integer linear program is an optimization problem of the form
\[ \mathop{\arg\min}\limits_{\textbf{x}} \left\{ \textbf{c}^T\textbf{x} \mid \textbf{Ax} \leq \textbf{b}, \textbf{l} \leq \textbf{x} \leq \textbf{u}, \textbf{x} \in \mathbb{Z}^\textit{p} \times \mathbb{R}^\textit{n-p} \right\} \]
where $\textbf{c} \in \mathbb{R}^\textit{n}$ is the objective coefficient vector, 
$\textbf{A} \in \mathbb{R}^{\textit{m} \times \textit{n}}$ is the constraint coefficient matrix, 
$\textbf{b} \in \mathbb{R}^\textit{m}$ is the constraint right-hand-side vector, 
$\textbf{l}, \textbf{u} \in \mathbb{R}^\textit{n}$ represent the lower and upper variable bound vectors respectively, 
and $\textit{p}$ is the number of integer variables.
The linear programming (LP) relaxation of a MILP is shown below
\[ \mathop{\arg\min}\limits_{\textbf{x}} \left\{ \textbf{c}^T\textbf{x} \mid \textbf{Ax} \leq \textbf{b}, \textbf{l} \leq \textbf{x} \leq \textbf{u}, \textbf{x} \in \mathbb{R}^\textit{n} \right\} \]
The LP solution provides a lower bound to the original MILP. Specifically, if the LP solution is subject to the integer constraint, then it is also a optimal feasible solution of the MILP. Otherwise, the LP relaxation is required to be decomposed into two sub-problems.
This is done by branching on a variable that does not obey the integrality constraint in the current LP solution.
The solving process terminates when the feasible regions cannot be decomposed anymore, and then a certificate of optimality or infeasibility can be provided respectively.

\subsection{Branching rules}

A key factor influencing the efficiency of B\&B algorithm is how to select a fractional variable to branch. 
In this part, three typical variable selection strategies are briefly described as follows.

The idea of Strong Branching is to evaluate which of the fractional candidate variables gives the best progress before actually branching on any of them.
For each candidate variable, this evaluation process is realized by solving the LP relaxations of the two sub-problems.
Thus, a huge amount of computation is required when adopting SB method.

Pseudocost branching is a sophisticated rule in the sense that it keeps a history of the success of the variables on which already has been branched~\cite{Benichou:3}.
$\Psi_j^+$ ($\Psi_j^-$) denotes the average unit objective gain taken over upwards (downwards) branching on $x_j$ in previous nodes.
Pseudocost branching at node $\textit{N}$ with LP relaxation solution $\check{x}$ consists in computing values:
\[ PC_j = score((\check{x}_j - \left \lfloor \check{x}_j \right \rfloor)\Psi_j^-,(\left \lceil \check{x}_j \right \rceil - \check{x}_j)\Psi_j^+) \]
and choosing the candidate variable with highest such value.
Compared with the SB method, PC method is simpler but faster.

Patel proposed an active constraint method that relies on estimating the impact pf the candidate variables on the active constraints in the current LP relaxation~\cite{patel2007active}.
This method aims to find the first feasible solution of MILP as quickly as possible.
The goal of the method is to select the branching candidate variable so that the LP relaxation optimum points for the two child nodes are as far apart as possible (in the sense of Euclidean distance).
This scheme can lead to significantly different solutions that one of the child nodes will be quite good while the other one is poor.
It is accomplished by choosing the candidate variable that most affects the active constraints at the parent node LP relaxation optimum.

%% file: content/solution.tex
\section{Overview of \Yordle}

\begin{figure*}[h!]
	\centering
	\includegraphics[height=8.4cm]{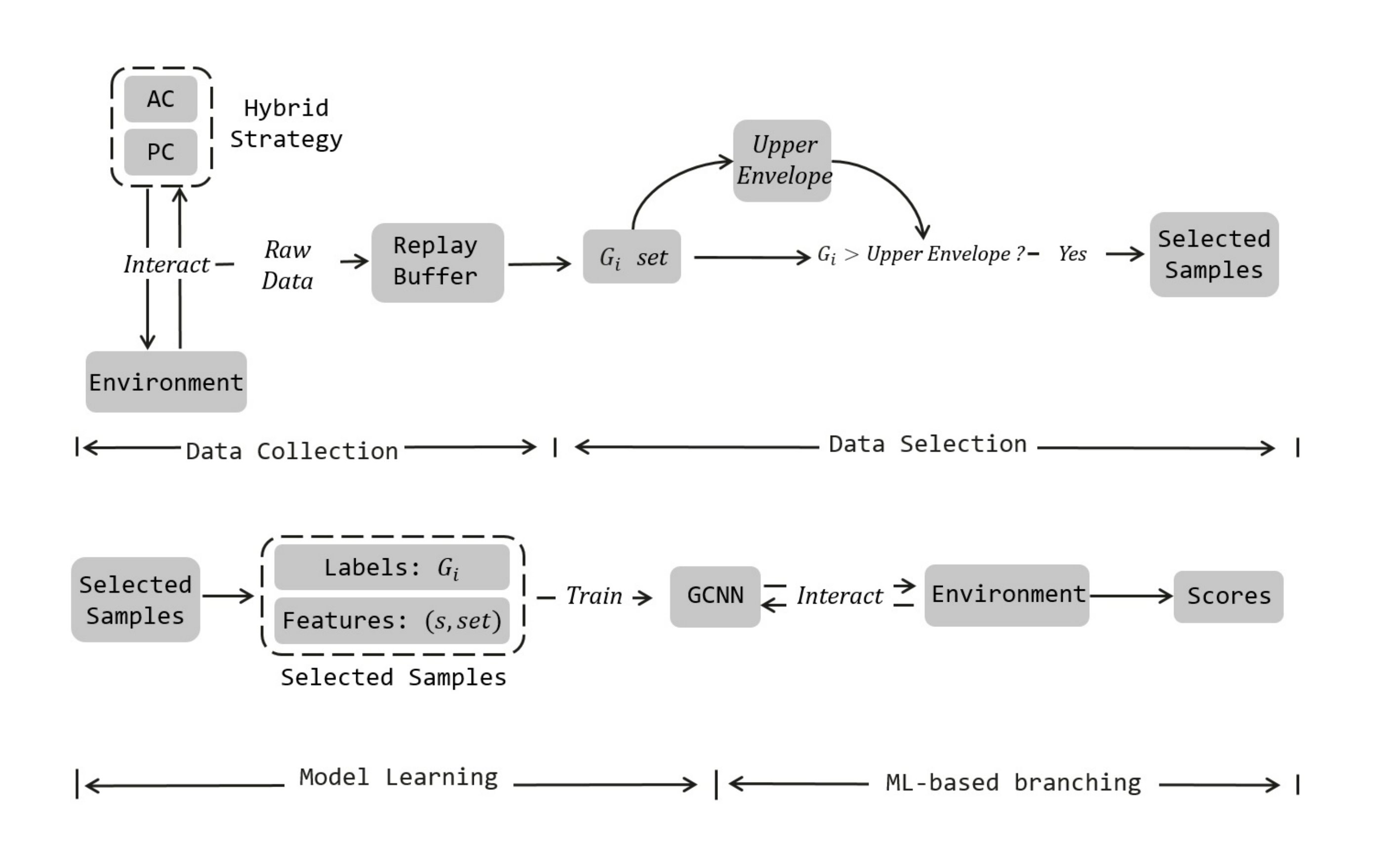}
	\caption{Framework of \Yordle. This framework consists of four phases, including data collection phase, data selection phase, model learning phase and ML-based branching phase. In the data collection phase, a hybrid strategy is adopted to collect expert data. And then in the data selection phase, we select the state-action pairs that possess higher cumulative reward. In the model learning phase, a graph convolutional neural network is adopted to train a branching strategy. Finally in the ML-based branching phase, the trained branching strategy is evaluated in parallel.}
	\label{fig: framework}
\end{figure*}
In our work, a framework, named \Yordle, is introduced for learning to branch based on BAIL. And it is proceeded in the following four phases:
\begin{enumerate}
    \item In the data collection phase, we collect the data with a hybrid strategy based on PC and AC strategies. Then the data is sent to the next phase.
    \item In the data selection phase, we develop a data selection method based on BAIL, and select the state-action pairs that possess higher cumulative reward. These selected state-action pairs are provided to the network model as expert demonstrations.
    \item In the model learning phase, similar to the Gasse's work, a graph convolutional neural network (GCNN) is adopted to exploit the natural bipartite graph representation of MILP problems. Differently, the criteria for generating candidate branching variables is cumulative reward, instead of instant reward.
    \item In the ML-based branching phase, a parallel evaluation platform is proposed to accelerate the instance evaluation process. Besides, inspired by the reinforcement learning, the evaluation criterion of the best strategy is chosen to be the cumulative reward rather than the cross-entropy loss.
\end{enumerate}
The framework is shown in Figure~\ref{fig: framework}. Next, we describe each of the phase in detail.

\subsection{Data collection}

In this phase, a hybrid strategy is adopted to collect the data to form a training dataset.
The setting of the hybrid strategy is as follows in order to increase the coverage of the dataset.
\[ strategy = \left\{ \begin{array}{cc}
PC & (db \leq DB_0, r \leq R_0) \ or \ (db > DB_0, r > R_0)\\
AC(r_1,r_2,r_3,r_4) & (db \leq DB_0, r > R_0) \ or \ (db > DB_0, r \leq R_0)
\end{array} \right. \]
where \textit{db} represents the value of dual bound of the instance, $ \textit{DB}_0 $ represents a certain value of dual bound, $ \textit{R}_0 \in (0,1) $ represents a sampling probability, $ r, r_1, r_2, r_3, r_4 \in [0,1] $ represent random quantities.

At each node $ N_i $, the training data comprises~\cite{prouvost2020ecole}:
\begin{itemize}
    \item Observation: A node bipartite graph representation of B\&B states used in Gasse et al., using the ecole.observation.NodeBipartite observation function. On one side of that bipartite graph, nodes represent the variables of the problem, with a vector encoding features of that variable. On the other side of the bipartite graph, nodes represent the constraints of the problem, similarly with a vector encoding features of that constraint. An edge links a variable and a constraint node if the variable participates in that constraint, that is, its coefficient is nonzero in that constraint. The constraint coefficient is attached as an attribute of the edge.
    \item Action set: The set of candidate variables.
    \item Action: The selected candidate variable to branch on.
    \item Reward: The reward is defined as the dual integral since the previous state, where the integral is computed with respect to the solving time.
    \item Next observation: The node bipartite graph representation of the next node.
    \item Next action set: The set of candidate variables corresponding to the next node.
    \item Done: The termination flag.
\end{itemize}
To demonstrate the framework of \Yordle in detail, those variables and sets are denoted by certain symbols, which are shown in Table~\ref{tab:symbol}:
\begin{table}[h]
    \centering
    \caption{Meaning of Symbols}
    \label{tab:symbol}
    \begin{tabular}{c c}
         \hline
         Symbols & Meaning \\
         \hline
         $obs$ & observation\\
         $set$ & action set\\
         $a$ & action\\
         $r$ & reward\\
         $obs^{\prime}$ & next observation \\
         $set^{\prime}$ & next action set\\
         $d$ & done \\
         \hline
    \end{tabular}
\end{table}

\subsection{Data selection}

In this phase, a data selection method based on BAIL algorithm~\cite{chen2019bail} is developed to obtain a better training dataset. The details are described as follows.

\begin{figure*}[h!]
	\centering
    \includegraphics[height=4.6cm]{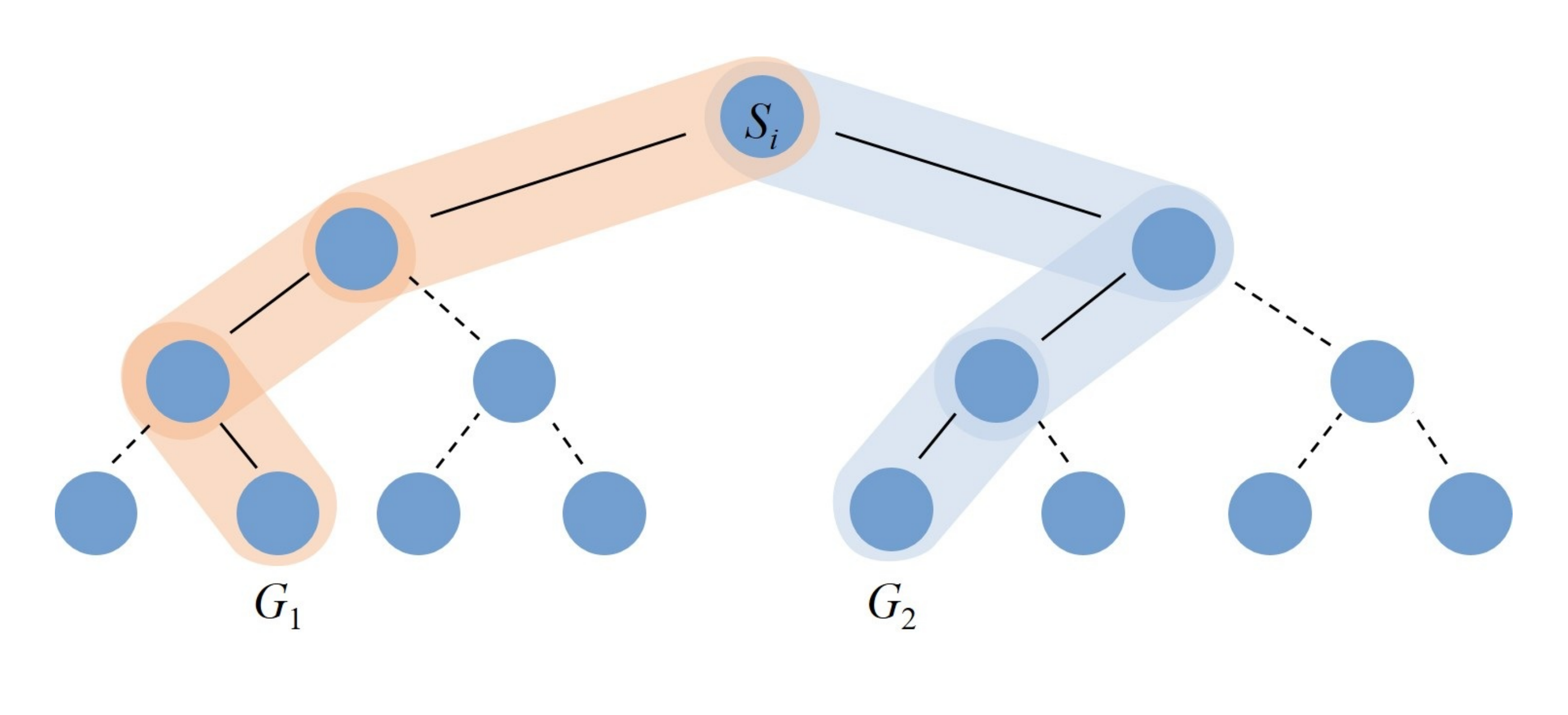}
    \caption{Cumulative Reward of the Forward Trace. All states $(obs,set)$ are firstly mapped to some low-dimensional states $s_i$. And then all the obtained low-dimensional states $s_i$ are stored into a buffer and traversed to get all possible paths. Then for each state, the cumulative rewards corresponding to all the paths starting from this state are calculated and finally the state-action pairs that possess higher cumulative reward are chosen for training.}
	\label{fig: reward}
\end{figure*}

A training dataset has been obtained in the data collection phase, which is denoted by $\mathcal{B} = \{ (obs_i, set_i,a_i,r_i,obs_i^{\prime},set_i^{\prime},d_i),i=1,...,m \}$, where subscript $i$ is a stamp to identify data points.
For each data point $i \in \{ 1,...,m \}$, we calculate the cumulative reward from the current state $(obs,set)$ to the end of the episode, which is shown as
\[ G_i = \sum_{t=i}^T{\gamma^{t-i}r_t} \]
where \textit{T} represents the step at which the episode ends for the episode that contains the $i_{th}$ data point, $\gamma$ represents the discount factor.
This process is demonstrated in Figure~\ref{fig: reward}.

In order to calculate the cumulative reward for a state $(obs,set)$, all the paths starting from $(obs,set)$ are required to be given. One way to do this is to put all states into a buffer and traverse them to get all possible paths. However, regarding the B\&B problem, the dimensionality of the state $(obs,set)$ is so large that it will take up too much storage resources. This is detrimental to our framework.
To address the problem, we introduce the following function that maps the state $(obs,set)$ to a low-dimensional vector $s$.
\[ f_{dim}(obs,set) = s \]
where $f_{dim}$ ensures that $(obs,set)$ and $s$ are in one-to-one correspondence.

For a fixed $\lambda \geq 0$, denote $\mathcal{G} = \{ (obs_i,set_i,G_i), i=1,...,m \}$.
Let $V_{\phi}(obs,set)$ denote a graph convolutional neural network characterized by $\phi = (w,b)$ that takes $obs_i$ and $set_i$ as input and outputs a real number to fit the cumulative reward.
Then $V_{\phi^\lambda}(obs,set)$ is regarded as a $\lambda$-regularized upper envelope for $\mathcal{G}$ if $\phi^\lambda$ is an optimal solution for the following constrained optimization problem:
\[ \min_\phi \sum_{i=1}^m {[V_\phi(obs_i,set_i) - G_i]^2 + \lambda \left \| w \right \|^2} \quad s.t. \quad V_\phi(obs_i,set_i) \geq G_i, i=1,...m \]

In the work of Chen et al.[11], an unconstrained optimization problem with a penalty loss function (with $\lambda$ fixed) is introduced to obtain an approximate upper envelope of the data $\mathcal{G}$, which is shown as
\[ L^K(\phi) = \sum_{i=1}^m {(V_\phi(obs_i,set_i)-G_i)^2 \{ \textbf{1}_{(V_\phi \geq G_i)} + K \cdot \textbf{1}_{(V_\phi < G_i)} \} + \lambda \left \| w \right \|^2} \]
where $K \gg 1$ represents the penalty coefficient, $\textbf{1}_{()}$ represents a indicator function.
For a certain finite \textit{K}, the penalty loss function will lead to an approximate upper envelope $V_\phi(obs_i,set_i)$.

Then we select all $(obs_i,set_i,a_i)$ pairs from the batch data set $\mathcal{B}$ such that
\[ G_i > x V_\phi(obs_i,set_i)\]
where \textit{x} is set such that the top $p\%$ of the data points are selected, where \textit{p} is a hyper-parameter.
In this paper, \textit{p} is set as 15.

\subsection{Model learning}

In this work, the imitation learning is applied to learn the policy from the selected data set $\mathcal{D}$.
Specifically, it is trained by minimizing the cross-entropy loss
\[ \mathcal{L}(\theta) = -\frac{1}{N} \sum_{(obs_i,set_i,a_i) \in \mathcal{D}} {\rm log} \pi_0(a \mid obs_i, set_i) \]

Then a graph convolutional neural network~\cite{gasse2019exact} is adopted to parametrize the candidate variable selection policy. 
In detail, the input of the GCNN model is the bipartite state representation $s_t = (\mathcal{G}, \textbf{C}, \textbf{V}, \textbf{E})$.
The graph convolution can be broken down into two successive passes, one from variables to constraints and one from constraints to variables, which are shown as
\begin{align*}
    c_i &\gets f_c(c_i, \sum_j^{(i,j) \in \varepsilon} {g_c(c_i,v_j,e_{i,j})}) \\
    v_i &\gets f_v(v_j, \sum_j^{(i,j) \in \varepsilon} {g_v(c_i,v_j,e_{i,j})})
\end{align*}
where $f_c$, $f_v$, $g_c$ and $g_v$ are two-layer perceptrons with ReLU activation functions.

Different with the work of Gasse et al.~\cite{gasse2019exact}, we produce the probability distribution over the candidate branching variables (i.e., the non-fixed LP variables) based on the cumulated rewards of each candidate variable, instead of the strong branching score. In this way, we can obtain a more effective and non-myopic policy.

\subsection{ML-based branching}
\begin{figure}[h!]
	\centering
	\includegraphics[height=4.9cm]{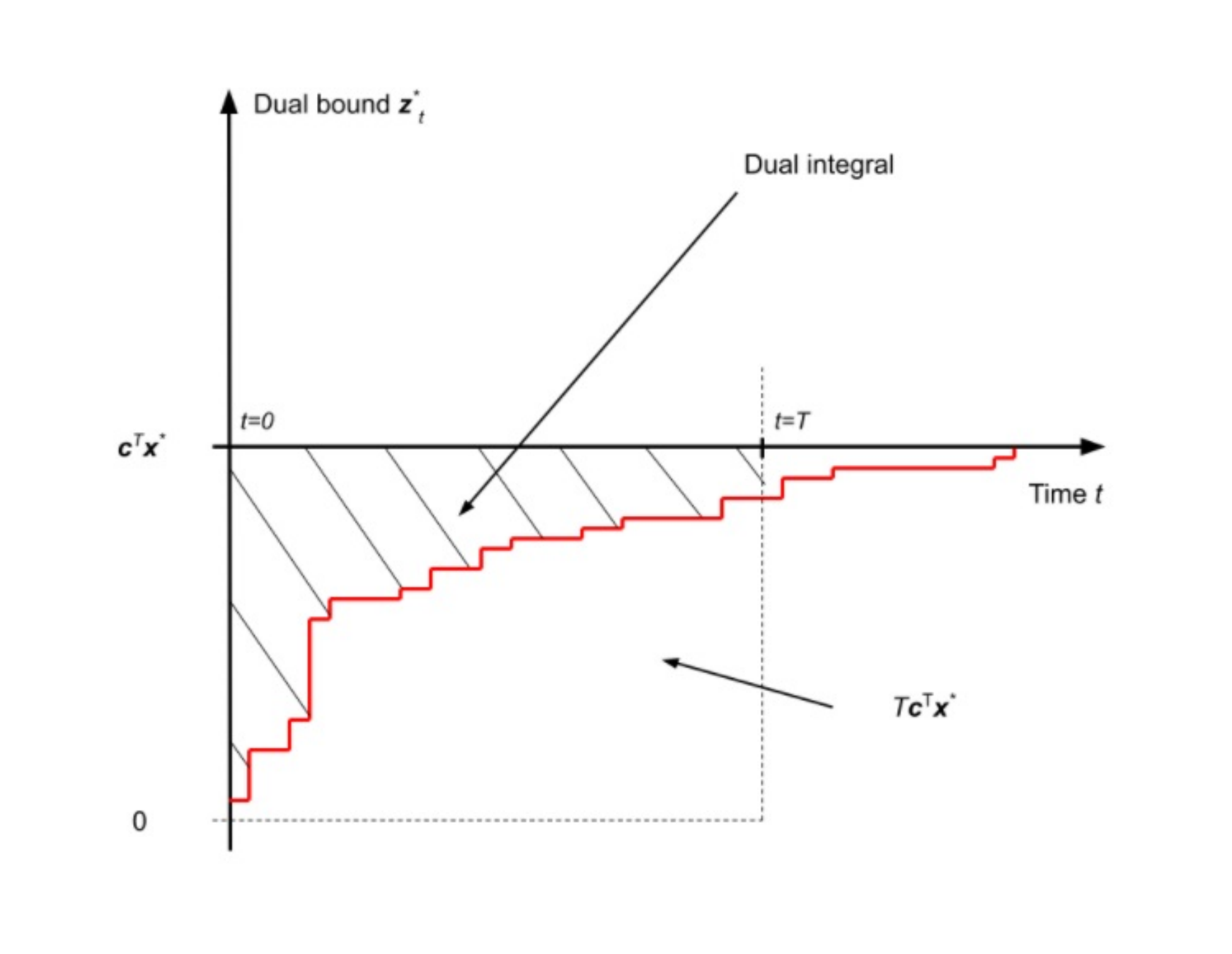}
	\caption{Dual Integral of a Instance}
	\label{fig: dual integral}
\end{figure}

In the evaluation stage, the trained GCNN model is used to replace the expert strategy to make decision while executing branch and bound. The environment will run a full-fledged branch-and-cut algorithm with SCIP, and the trained model will only control the solver's branching decisions. Also, all primal heuristics will be deactivated, so that the focus is only on proving optimality via branching. 
In each instance of MILP problem, the metric to evaluate a branch and bound algorithm is dual integral, which is expressed as follow:
\begin{align*}
    T\mathbf{c}^{T}\mathbf{x^{*}}-\int_{t=0}^{T}\mathbf{z^{*}}dt
\end{align*}
where $\mathbf{z^{*}}$ is the best dual bound at time t, and $ T\mathbf{c}^{T}\mathbf{x^{*}}$ is an instance-specific constant that depends on the optimal solution $\mathbf{c}^{T}\mathbf{x^{*}}$. The dual integral is to be minimized, and takes an optimal value of 0. Equivalently, the cumulative reward, which is given by a constant minus the dual integral, is to be maximized in our task. Intuitively, Figure~\ref{fig: dual integral} shows the dual integral.
Considering that the valid set includes a large number of instances, a parallel evaluation mechanism is introduced to accelerate the evaluation process. However, despite the shortened evaluation time, it also misleads our experimental results to some extent, as will be explained in the experiment.

%% file: content/evaluation.tex
\section{Evaluation}

\subsection{Competition results}

In the competition, the performance of submitted models are evaluated in three problem benchmarks from diverse application areas, namely `Item placement', `Load Balancing' and `Anonymous' problem. Each of the benchmark dataset consists of many MILP instances, which are split into two distinct collections: train and valid. The cumulative reward of each task is given by performing the model in a hidden test datasets. We submitted \Yordle~ and named our team as qqy.

Due to the time constraint of this competition, we only implement \Yordle~ over the dataset of `Item Placement' and `Anonymous'. It takes a lot of time to collect data on the dataset `Load Balancing', and it is too late for us to implement our framework on the dataset `Load Balancing' because designing this framework consumed a lot of our time. In fact, we adopt the method proposed by Gasse et al.~\cite{gasse2019exact} on the dataset `Load Balancing'. Therefore, in this part, we only care about the competition results on the `Item Placement' and `Anonymous'. And some discussions about the `Load Balancing' will be given in the 
following.

Table~\ref{tab:res_on_ip_and_ap} shows the scores of the top 10 leading teams in global leaderboard for the `Item Placement' and `Anonymous', alongside with the ranks of them.
Table~\ref{tab:res_on_ip_and_ap_st} shows the scores of the top 5 leading teams in student leaderboard for the `Item Placement' and `Anonymous', alongside with the ranks of them.
It is demonstrated from the result that \Yordle~ is effective on the `Item Placement' and `Anonymous' problems.

\begin{table*}[t!]
    \centering
    \caption{Results on the `Item Placement' and `Anonymous' Problems in Global Leaderboard}
    \label{tab:res_on_ip_and_ap}
    \begin{tabular}{c c c|c c c}
    \hline
    \quad & item placement & \quad & \quad & anonymous & \quad \\
    \hline
    team name & cum. reward & rank & team name & cum. reward & rank \\
    \hline
    Nuri & 6684.00 & 1 & Nuri & 27810782.42 & 1 \\
    EI-OROAS & 6670.30 & 2 & \textbf{qqy} & \textbf{27221499.03} & \textbf{2} \\
    EFPP & 6487.53 & 3 & null\_ & 27184089.51 & 3 \\
    lxj24 & 6443.55 & 4 & EI-OROAS & 27158442.74 & 4 \\
    ark & 6419.91 & 5 & DaShun & 27151426.15 & 5 \\
    \textbf{qqy} & \textbf{6377.23} & \textbf{6} & KEP-UNIST & 27085394.46 & 6 \\
    KAIST\_OSI & 6196.56 & 7 & lxj24 & 27052321.48 & 7 \\
    nf-lzg & 6077.72 & 8 & THUML-RL & 26824014.00 & 8 \\
    Superfly & 6024.20 & 9 & KAIST\_OSI & 26626410.86 & 9 \\
    Monkey & 5978.65 & 10 & Superfly & 26373350.99 & 10 \\
    \hline
    \end{tabular}
\end{table*}

\begin{table*}[t!]
    \centering
    \caption{Results on the `Item Placement' and `Anonymous' Problems in Student Leaderboard}
    \label{tab:res_on_ip_and_ap_st}
    \begin{tabular}{c c c|c c c}
    \hline
    \quad & item placement & \quad & \quad & anonymous & \quad \\
    \hline
    team name & cum. reward & rank & team name & cum. reward & rank \\
    \hline
    lxj24 & 6443.55 & 1 & \textbf{qqy} & \textbf{27221499.03} & \textbf{1} \\
    ark & 6419.91 & 2 & null\_ & 27184089.51 & 2 \\
    \textbf{qqy} & \textbf{6377.23} & \textbf{3} & lxj24 & 27052321.48 & 3 \\
    KAIST\_OSI & 6196.56 & 4 & THUML-RL & 26824014.00 & 4 \\
    nf-lzg & 6077.72 & 5 & KAIST\_OSI & 26626410.86 & 5 \\
    \hline
    \end{tabular}
\end{table*}

\begin{figure*}[h!]
	\centering
	\subfigure[Item Placement]{\includegraphics[width=6.6cm]{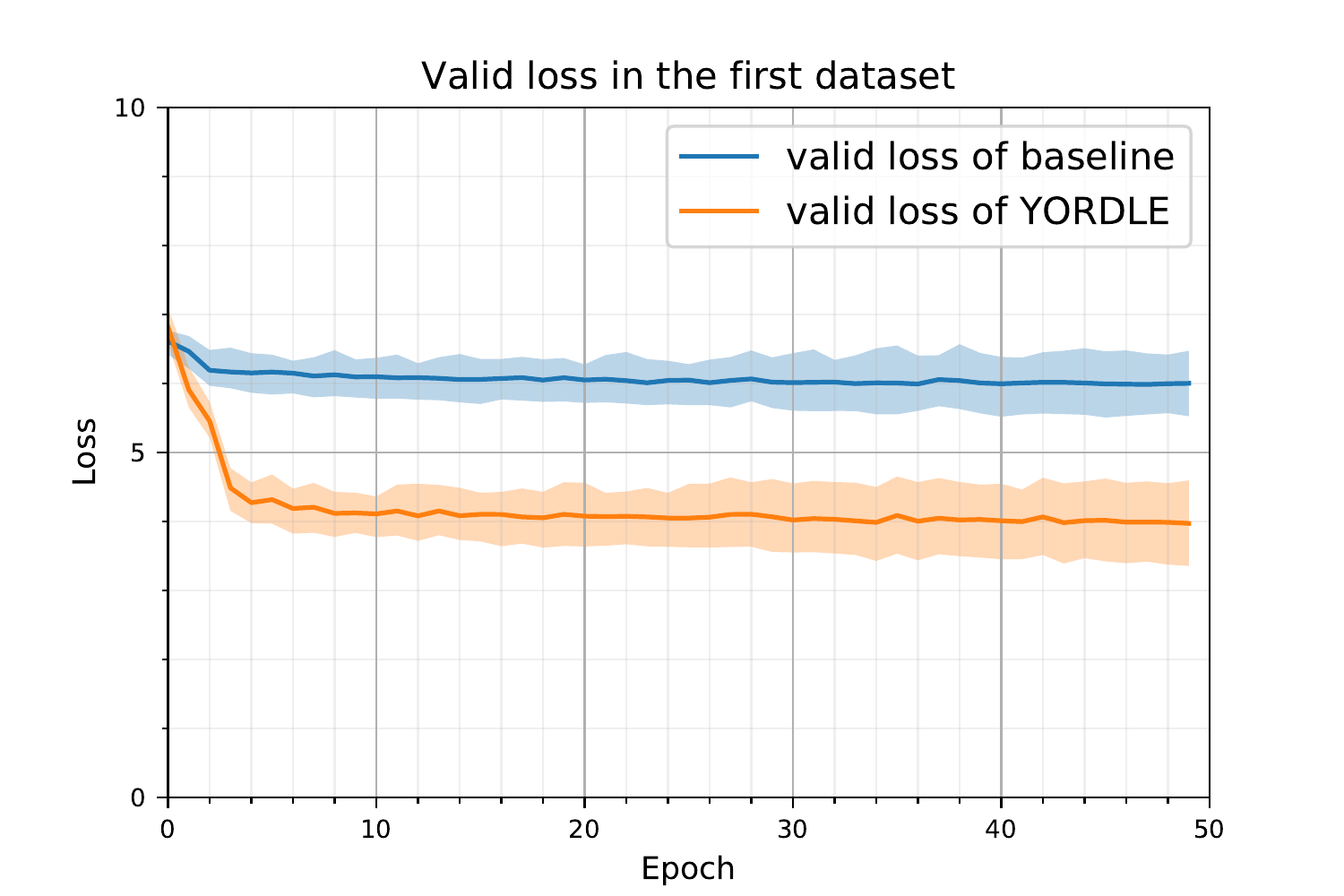}}
	\subfigure[Anonymous]{\includegraphics[width=6.6cm]{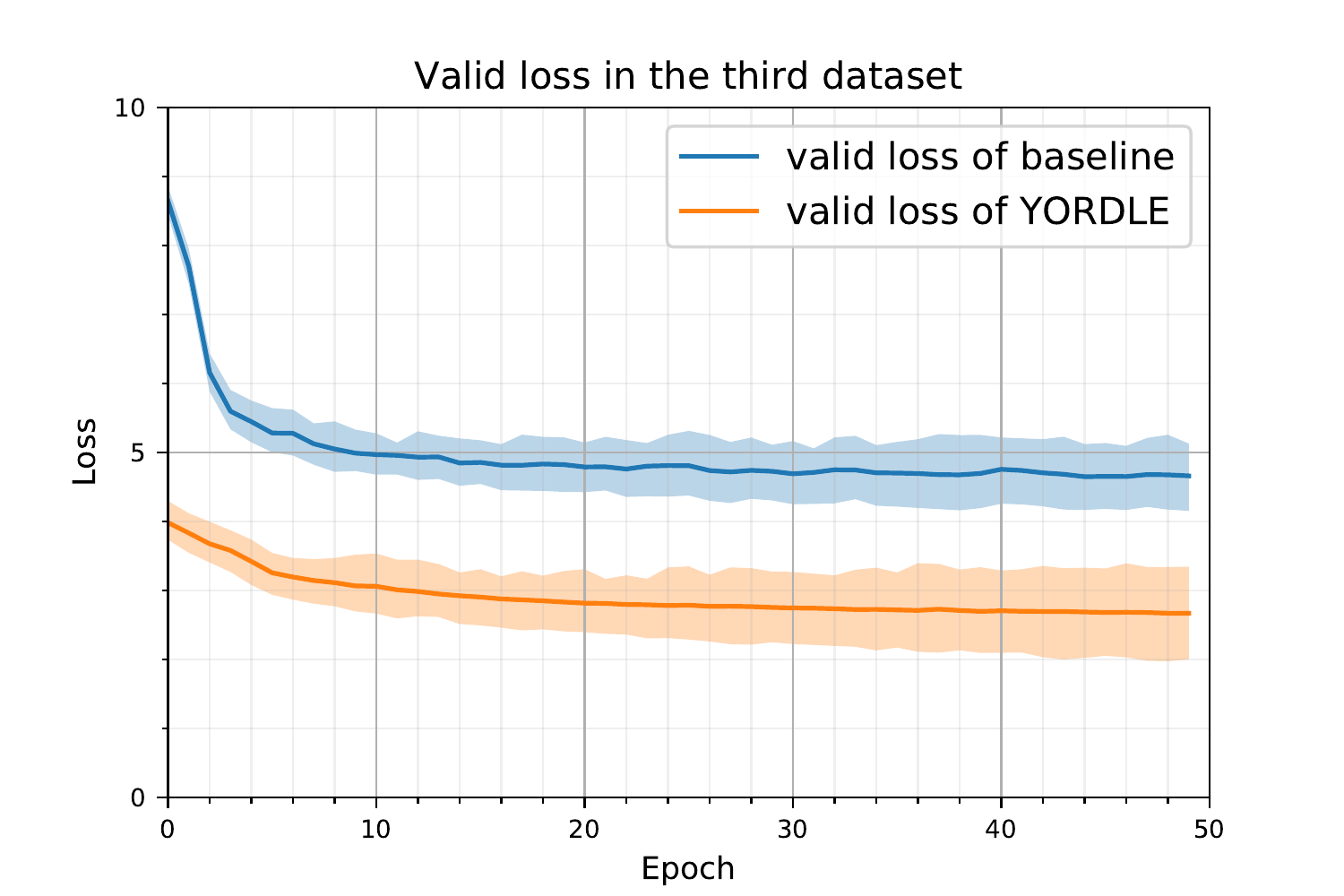}}
	\caption{Valid Loss on the Item Placement and Anonymous Problems}
	\label{fig: loss_over_IP_and_AP}
\end{figure*}

Note that compared with the imitation learning method proposed by Gasse et al.~\cite{gasse2019exact}, it can be observed in Figure~\ref{fig: loss_over_IP_and_AP} that our method converges faster and possesses much smaller valid loss, which makes it easier to train. Therefore, we believe that it is a small but efficient imitation learning framework for learning to branch.

\subsection{Supplementary Experiment}

After the competition, we completed the training of \Yordle~ on the dataset `Load Balancing'. Considering that our device may be different with that of the competition, we evaluate the cumulative reward of both our submitted strategy (i.e. qqy) and the one trained by means of \Yordle~ on the local device, which are shown in Table~\ref{tab:eval_on_local_device}. We consider the following assumptions to be reasonable.
\begin{align*}
    \frac{R_{ldb}}{R_{loc}} = CONSTANT
\end{align*}
where $ R_{ldb} $ and $ R_{loc} $ represent the cumulative reward of a strategy on the leaderboard and the local device respectively. In this way, it can be inferred that the cumulative reward of \Yordle~ on the leaderboard is about 630793.45. In fact, the strategy trained by means of \Yordle~ can be ranked $ 6_{th} $ in the global leaderboard and ranked $ 2_{nd} $ in the student leaderboard.
Table~\ref{tab:res_on_lbp_gl} shows the scores of the top 10 leading teams in global leaderboard for the `Load Balancing' problem, alongside with the ranks of them.
Table~\ref{tab:res_on_lbp_sl} shows the scores of the top 5 leading teams in student leaderboard for the `Load Balancing' problem, alongside with the ranks of them.
It is demonstrated from the result that \Yordle~ is effective on the dataset `Load Balancing' as well.

\begin{table*}[h!]
    \centering
    \caption{Evaluation Results of qqy and \Yordle~ on the Local Device}
    \label{tab:eval_on_local_device}
    \begin{tabular}{c c c}
    \hline
    official cum. reward of qqy & local cum. reward of qqy & local cum. reward of \Yordle~ \\
    \hline
    630557.31 & 630629.31 & 630865.48 \\
    \hline
    \end{tabular}
\end{table*}

\begin{table}[h!]
    \centering
    \caption{Results on the `Load Balancing' Problem in Global Leaderboard}
    \label{tab:res_on_lbp_gl}
    \begin{tabular}{c c c}
    \hline
    team name & cum. reward & rank \\
    \hline
    EI-OROAS & 631744.31 & 1 \\
    KAIST\_OSI & 631410.58 & 2 \\
    EFPP & 631365.02 & 3 \\
    DaShun & 630898.25 & 4 \\
    blueterrier & 630826.33 & 5 \\
    \textbf{\Yordle~} & \textbf{630793.45} & \textbf{6} \\
    Nuri & 630787.18 & 7 \\
    gentlemenML4CO & 630752.94 & 8 \\
    comeon & 630750.66 & 9 \\
    Superfly & 630746.96 & 10 \\
    \hline
    \end{tabular}
\end{table}

\begin{table}[h!]
    \centering
    \caption{Results on the `Load Balancing' Problem in Student Leaderboard}
    \label{tab:res_on_lbp_sl}
    \begin{tabular}{c c c}
    \hline
    team name & cum. reward & rank \\
    \hline
    KAIST\_OSI & 631410.58 & 1 \\
    \textbf{\Yordle~} & \textbf{630793.45} & \textbf{2} \\
    comeon & 630750.66 & 3 \\
    qqy & 630557.31 & 4 \\
    ark & 630414.45 & 5 \\
    \hline
    \end{tabular}
\end{table}

\begin{figure}[!h]
	\centering
	\includegraphics[width=6.6cm]{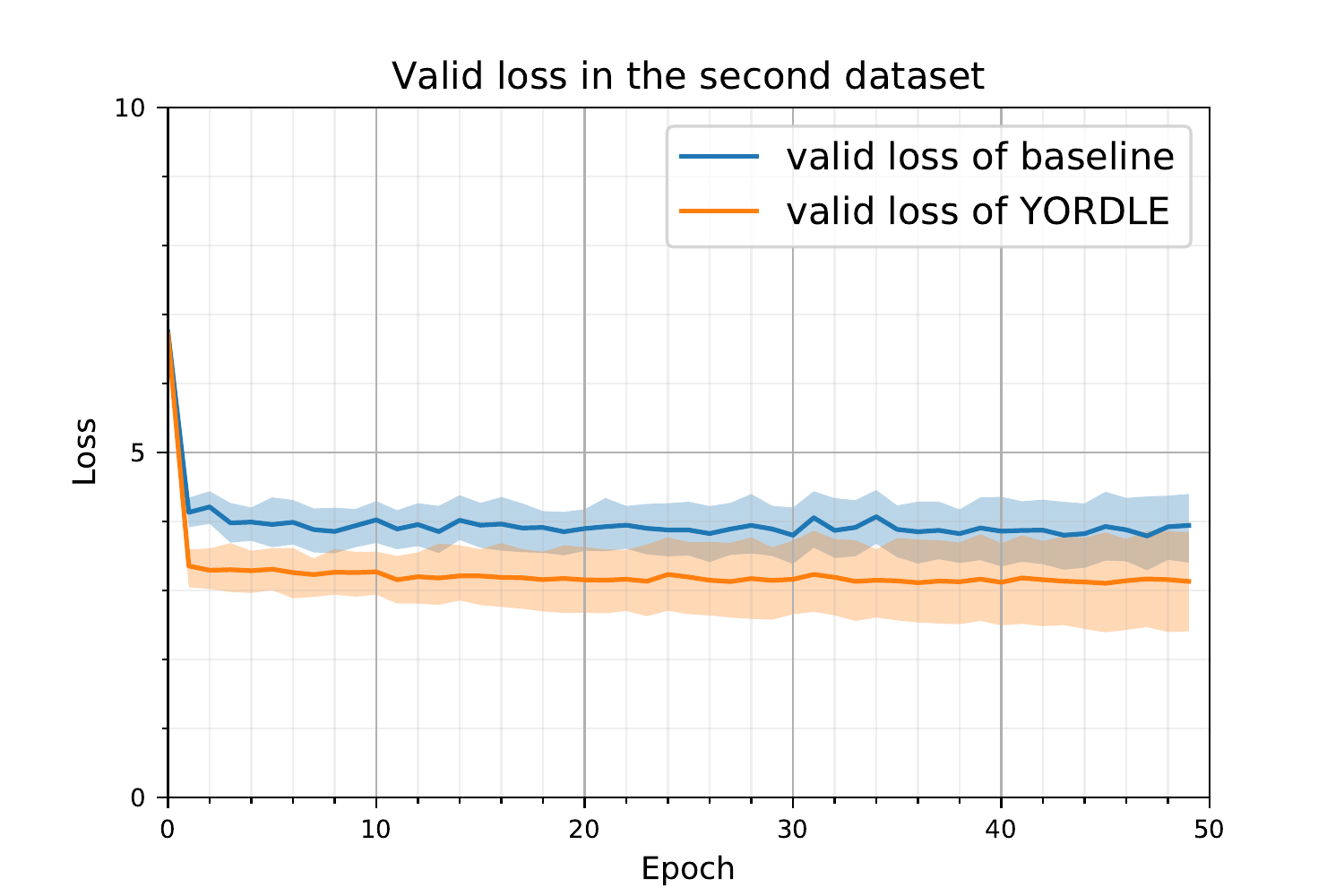}
	\caption{Valid Loss on the `Load Balancing' Problems}
	\label{fig: loss_lbp}
\end{figure}

Also, the valid loss is given in Figure~\ref{fig: loss_lbp}. It is still observed that our method converges faster and possesses smaller valid loss on the `Load Balancing' problem.

%% file: content/conclusion.tex
\section{Conclusion}
In this paper, we introduce our framework \Yordle, which achieved the most excellent result in student leaderboard of the ML4CO competition. Generally, \Yordle~ adopts an imitation learning method, which leverages a hybrid expert strategy to collect data and selects state-action pairs with higher cumulative reward to train a GCNN as the branch and bound policy. \Yordle~ is flexible to be extended and easy to implement while requiring shorter time and fewer data to train. During the competition, we also have several observations. First, the existing heuristic method solving MILP problem may be not the best solution. Especially, the strong branch rule, which is commonly believed to be the `golden standard', is inferior to the sampling rule we used in this competition. Last, one of the bottlenecks to implement a machine leaning method on branch and bound is the representation capability of the model, especially when the dataset is large. This conclusion is drawn from the observation that ML methods might not outperform the random algorithm when learning from the large dataset. Normally, extending the size of the model helps to improve the representative ability. However, it also leads to a slower solving process. Thus, we have to deal with a trade off between efficiency and effectiveness.